\documentclass{psp}
\usepackage[english]{babel}
\usepackage[T1]{fontenc}
\usepackage[latin1]{inputenc}
\usepackage{graphicx}
\usepackage{amsmath,amssymb}
\usepackage{url}

 \newtheorem{thm}{Theorem}
 \newtheorem{cor}[thm]{Corollary} 
 \newtheorem{lem}[thm]{Lemma} 
 \newtheorem{defn} [thm] {Definition}
 
 \newtheorem{conj} [thm]{Conjecture}
 \newtheorem{prop}[thm]{Proposition} 
 
 \newtheorem{rem*}{Remark}

\begin{document}

\title[Polynomials of Annular Surfaces and their Boundary Links]
{On Knot Polynomials of Annular Surfaces and their Boundary Links}
\makeatletter
\author[Hermann Gruber]{Hermann Gruber\thanks{Most of the work was done while the author was at Institut f\"ur Informatik,
Ludwig-Maximilians-Universit\"at M\"unchen,
Oettingenstra\ss{}e~67, D-80538 M\"unchen, Germany.}\\
Institut f\"ur Informatik, 
Justus-Liebig-Universit\"at Giessen\addressbreak
Arndtstr. 2, 35392 Giessen, Germany\addressbreak
e-mail: \textup{\texttt{hermann.k.gruber@informatik.uni-giessen.de}}
}
\receivedline{Received \textup{7} March \textup{2008;}
              revised \textup{10} October \textup{2008}}
\makeatother
\maketitle
\begin{abstract}
Stoimenow and Kidwell asked the following question:
Let $K$ be a non-trivial knot, and let $W(K)$ be a 
Whitehead double of $K$. Let $F(a,z)$ be the 
Kauffman polynomial and $P(v,z)$ the skein polynomial.
Is then always $\max\deg_z P_{W(K)} - 1 = 2 \max\deg_z F_K$?
Here this question is rephrased in more general terms 
as a conjectured relation between the maximum $z$-degrees 
of the Kauffman polynomial of an annular
surface $A$ on the one hand, and the Rudolph polynomial on 
the other hand, the latter being defined as a certain
M\"obius transform of the skein 
polynomial of the boundary link $\partial A$.
That relation is shown to hold for 
algebraic alternating links, thus simultaneously solving  
the conjecture by Kidwell and Stoimenow and a related 
conjecture by Tripp for this class of links.
Also, in spite of the heavyweight definition of the 
Rudolph polynomial $\{K\}$ of a link $K$, 
the remarkably simple formula 
$\{\bigcirc\}\{L\#M\}=\{L\}\{M\}$
for link composition is established. 
This last result can be used to reduce the conjecture 
in question to the case of prime links. 
\end{abstract}

\section{Preliminaries.}

We assume the reader is familiar with basic notions in knot theory
as contained, e.g., in~\cite{BZ85}. To fix the notation, we recall
some notions on framed links and framed link polynomials
that are, in part, not yet contained in~\cite{BZ85}.

\subsection{Framed Links.}\label{framed-links}

Let $A$ be an orientable 
annulus $S^1 \times I$ 
contained in an unknotted solid torus $V$ in $S^3$ 
such that $A$ is parallel to the core of $V$, 
and let $N(K)$ denote the tubular neighborhood
of a knot $K$ in $S^3$. Let $h: V \to N (K)$ be an orientation 
preserving homeomorphism taking the meridian disk of $V$ to the 
meridian disk of $N(K)$ and the core of $V$ to $K$.
The image $h(A)$, a knotted annulus, is called 
the {\em framed knot} $A(K,f)$ whose framing $f\in \mathbb{Z}$ is the 
linking number of $h(\ell)$ and $K$, with $\ell$ being the preferred 
longitude of $V$. The mapping $h$ naturally induces an orientation
on $h(A)$ and its boundary $\partial h(A)$. The notation $A(K,f)$
is justified by the fact that, up to ambient isotopy, $A(K,f)$
depends only on $K$ and $f$. A change of orientation 
yields $A(K,f)=A(-K,f)=-A(K,f)$. 
The boundary $\partial A(K,f)$ is an oriented $2$-component link, 
sometimes referred to as the {\em $f$-twisted double} of $K$. 

The above notions readily generalize to links $L$ in $S^3$ with $\mu$
components $C_1,\ldots C_\mu$, by specifying corresponding homeomorphisms 
$h_i: V_i \to N(C_i)$, with unknotted solid tori $V_i$ as preimages. 
A {\em framed link} is an annular surface $A(L,f)$ specified by a 
pair $(L,f)$ as above, but now the framing $f$ maps $L$ to a vector 
in $\mathbb{Z}^\mu$ specifying the framing for each component $C_i$.
A {\em sublink} $L'$ of $L$ is obtained from $L$ by erasing zero or more
components from $L$; A framing $f$ of $L$ naturally induces a framing
on the sublinks. As long as there is no risk of confusion, this framing will be also 
denoted by $f$.   

For a link $L$, let $D = D(L)$ be a diagram of $L$ in $\mathbb{R}^2$.
We call a crossing in $D$ {\em homogeneous} if it involves only 
one component, and {\em heterogeneous} in the other case. 
We say that the crossing
$\vcenter{\hbox{\includegraphics{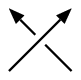}}}$  has positive sign,
and the crossing $\vcenter{\hbox{\includegraphics{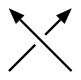}}}$ has negative sign.
For a link diagram $D(L)$, the {\em writhe} $w(D)$ is the overall sum of the crossings in $D$,
the {\em total linking number} $t(L)$ is half the sum of signs of the heterogeneous crossings
in $D$, and the {\em self-writhe} $sw(D)$ is the sum of signs of the homogeneous crossings in $D$.
The total linking number is invariant under ambient isotopy, and writhe and self-writhe
are regular isotopy invariants. Each link diagram $D = D(L)$ naturally induces
a framing~$f$ onto the surface $A(L,f)$. For a component $C$ of $L$, we have 
$f(C) = sw(D|_{C})$, where $D|_C$ is the diagram of~$C$ obtained by erasing
all other components from~$D$. 
The self-writhe is sometimes also called {\em total framing}~\cite{Rud90}.

\subsection{Link polynomials.}

The {\em skein polynomial} $P(v,z)$, defined first in~\cite{HOMFLY,PT}
is a Laurent polynomial in two variables $v$ and $z$. It is an ambient isotopy 
invariant of oriented knots and links and is defined {\em via} link diagrams
by  the (skein) relation 
 $$v^{-1}\vcenter{\hbox{\includegraphics{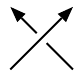}}}
-v\vcenter{\hbox{\includegraphics{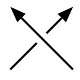}}}=z\vcenter{\hbox{\includegraphics{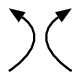}}},$$
and agreeing that the polynomial is equal to $\frac{v^{-1}-v}{z}$ 
for the unknot, the latter being denoted by $\delta_P$ as a shorthand. 
In addition, we let the polynomial be equal to $1$ on the empty diagram. 
It is not hard to see that then taking the split union of link 
diagrams amounts to multiplying their skein polynomials, also if 
one of the diagrams is empty.

Similarly, the {\em framed Kauffman polynomial} $L_{K,f}(a,x)$, defined 
in~\cite{Kau90}, is an ambient isotopy invariant for framed 
links $A(L,f)$ defined via the following relations, and
the convention that the polynomial equals $\frac{a^{-1}+a}{x}-1$ 
for the unknotted annulus with framing $0$ and $1$ for the empty diagram:

$$\vcenter{\hbox{\includegraphics{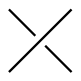}}}+\vcenter{\hbox{\includegraphics{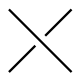}}}=
x\left(\vcenter{\hbox{\includegraphics{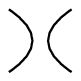}}}+\vcenter{\hbox{\includegraphics{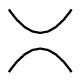}}}\right);$$
$$\vcenter{\hbox{\includegraphics{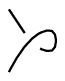}}}=a~\vcenter{\hbox{\includegraphics{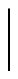}}};$$
the mirror image of this kink can be resolved by multiplying with $a^{-1}$.

The above diagrams are to be interpreted as parts of band diagrams 
of an annular surface that appears untwisted, 
and are denoted by  $K_+$, $K_-$, $K_\infty$ and $K_0$ respectively. 
Similar to above, write $\delta_L=\frac{a^{-1}+a}{x}-1$
for the polynomial of an unknotted annulus with framing $0$.
The framed Kauffman polynomial can be normalized to an 
unoriented link invariant, by letting $U_K =a^{-sw(K)}L_K$, or
to a semi-oriented link invariant, by letting $F_K =a^{-w(K)}L_K$ \cite{Kau90}. 

\begin{rem*}
Note that the normalization used here differs from the convention in 
some standard works on link polynomials, where these polynomials
equal $1$ on a simple closed circle, resp. an unknotted annulus 
with framing $0$, rather than on the empty diagram.
The author deliberately sticks to this convention because many of the formulas
appearing in this paper would become otherwise much more cumbersome. 
\end{rem*}

Yamada defined the following linear operator on polynomial 
invariants for framed links~\cite{Yam89}:
Let $R_K$ be a polynomial invariant for framed links $K$ under ambient isotopy, 
and let~$\lambda$ be a scalar. Let $|K|$ denote the number of components of~$K$.
Then the invariant~$(R_K\oplus\lambda)$ is defined 
as $\sum_{K'}\lambda^{|K|-|K'|}R_{K'}$,
where the summation index $K'$ runs over all framed sublinks of $K$,
including the empty link.

Now the {\em framed Rudolph polynomial} $\{K,f\}$ is defined as
$\{K,f\} = (P_{\partial A(K,f)}\oplus(-1))$.
Using~$[K,f]$ as a shorthand for the skein polynomial~$P_{\partial A(K,f)}$,
this is the M\"obius transform $\{K,f\}=\sum_{K'}(-1)^{|K|-|K'|}[K',f]$, 
as originally defined in~\cite{Rud90}. Note that for the special case of 
knots, we have $\{L,f\}=[L,f]-1$.
Write $A(L,0)$ if the framing for each annulus 
has the value~$0$, let $[L]:=[L,0]$ and $\{L\}=\{L,0\}$.
In~\cite[prop. 2(5)]{Rud90}, it is noted
that $\{L,f\}=v^{2sw(L)}\{L\}$, and $\{L\}$ is an ambient isotopy
invariant for the link $L$, that is, independent of the framing~$f$.

\subsection{Inequalities for Link Polynomials}

We briefly recall two known inequalities for the skein polynomial and
the Kauffman polynomial. Morton proved the following inequality:
\begin{thm}[Morton's Inequality]
Let $K$ be a link. For any link diagram $D$ of $K$, let $c(D)$ denote the number 
of crossings and $s(D)$ denote the number of Seifert circles obtained from
smoothing out all crossings, i.e. replacing each 
$\vcenter{\hbox{\includegraphics{fig/positive}}}$ and 
$\vcenter{\hbox{\includegraphics{fig/negative}}}$ with 
$\vcenter{\hbox{\includegraphics{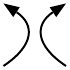}}}$.
Then $\max\deg_z P_K \le c(D) - s(D).$
\end{thm}
Kidwell proved an inequality relating properties of link diagrams 
to the maximum $z$-degree of a specialization of the (at that time yet
undiscovered) Kauffman polynomial~\cite{Kid87};
Thistlethwaite later showed that essentially the same proof 
applies to the more general Kauffman polynomial~\cite{Thi88}:
\begin{thm}[Kidwell's Inequality]\label{thm:Kidwell}
Let $K$ be a link. For any link diagram $D$ of $K$, let $c(D)$ denote the number 
of crossings and $b(D)$ denote the maximal bridge length, i.e. the longest 
sequence of consecutive crossing overpasses in $D$. Then
$\max\deg_x F_K \le c(D) - b(D)$. Moreover, equality holds if $D$ is a reduced 
alternating and prime diagram.
\end{thm}

\section{Reformulations of Rudolph's Congruence Theorem.}

Rudolph stated a curious theorem \cite{Rud90}, 
which in part motivates the present study, 
relating the Kauffman and Rudolph polynomials with coefficients reduced modulo~$2$.
We start this investigation by taking a closer look at his
Congruence Theorem:

\begin{thm}[Congruence Theorem]
Let $K$ be any link. Then
$$F_K(v^{-2},z^2)\equiv_2 v^{4t(K)}\{K\}.$$
\end{thm} 

We can reformulate the theorem in terms of the unoriented
and the framed Kauffman polynomial, rendering the
congruence more natural. In particular, the congruence can also 
be used to relate the framed versions of these invariants.
We also present versions involving a M\"obius transform of the 
Kauffman polynomial, and the polynomial $[L,f]$, which unveil
a certain symmetry in Rudolph's congruence.

\begin{prop}[Unoriented, framed, M\"obius transformed and Dubrovnik version]\label{prop:congruence-transform}
The following are equivalent reformulations of the congruence theorem:
\begin{eqnarray*} 
U_K(v^{-2},z^2)& \equiv_2& \{K\} \\
L_{(K,f)}(v^{-2},z^2)& \equiv_2& \{K,f\} \\
~[K,f] & \equiv_2 & ({L}_{(K,f)}\oplus \pm 1)(v^{-2},z^2)  \\
~[K,f] & \equiv_2 & (D_{(K,f)}\oplus \pm 1)(v^{-2},-z^2) 
\end{eqnarray*}
where the sum in the second and the third congruence runs over all
non-empty sublinks $K'$ of $K$.
\end{prop}
\begin{proof}
First, the signs in the $\pm 1$ terms in the congruences 
are immaterial since we work with coefficients reduced modulo~$2$.
For every framed link $(K,f)$, we have $w(K)=sw(K)+2t(K)$. Since
$F_K(a,x)=a^{-(sw(K)-2t(K))}L_K =a^{-2t(K)}U_K$,
we see that $F_K(v^{-2},z^2)=v^{4t(K)}U_K(v^{-2},z^2)$,
so that we can rewrite the congruence theorem as in the 
first congruence. For the second, we have
$F_K(v^{-2},z^2)=v^{4t(K)}v^{2sw(K)}L_{(K,f)}(v^{-2},z^2)$.
From~\cite[Prop.~2(5)]{Rud90},
we deduce $\{K,f\}=v^{-2sw(K)}\{K\}$. Thus, the second 
congruence is also equivalent to the original one.
For the third congruence, we use the fact that
$([K,f]\oplus0) = ([K,f]\oplus(-1)\oplus1) = \sum_{K'}\{K',f\}$, 
where the sum runs over all sublinks 
of~$K$~\cite[Prop.~1]{Yam89}; finally the 
fourth congruence is obtained by 
the Interchange Formula (see \cite{Lic88})
$$L_{(K',f)}(v^{-2},z^2)= (-1)^{w(K',f)+|K'|}D_{(K',f)}(v^{-2},-z^2)$$
for switching between the Kauffman and the Dubrovnik polynomial.
\end{proof}

It can be readily seen that the congruence theorem would equally hold
when using the variant $([K,f]\oplus1)$ instead 
of $\{K,f\} = ([K,f]\oplus(-1))$. 
As we will see later, however, the 
current choice is rewarded by an eminently simple formula for the Rudolph 
polynomial of composite links. 
We note on the fourth version of the congruence that
Yamada  proved in~\cite{Yam89} the remarkable one-variable {\em equality}
$$[L,f](t,t^{\frac{1}{2}}-t^{-\frac{1}{2}}) = (D_{(K,f)}\oplus 1)(t^{-2},t^{-1}-t),$$
The left-hand side of the above equation is the 
Jones polynomial of $\partial A(K,f)$ from \cite{Jon85}, 
normalized to value~$1$ at the empty diagram.

\section{Problem Statement and Preliminary Results}\label{sec:statement}

Kidwell and Stoimenow posed the following question
relating the $x$-degree of the Kauffman polynomial 
of a knot and the $z$-degree of its Whitehead 
double~(\cite{KidSto03}, see also~\cite{Oht04}):

\begin{conj}[Kidwell, Stoimenow]\label{conj:Kidwell}
Let $K$ be a non-trivial knot, and let $W(K)$ be a 
Whitehead double of $K$.
Is then always $\max\deg_z P_{W(K)} - 1 = 2 \max\deg_x F_K $?
\end{conj}

They note that their conjecture holds for all prime knots with up to $11$
crossings, by exhaustive calculation. A related conjecture
by Tripp states that the canonical genus of a Whitehead double
of every knot $K$ coincides with the minimal crossing number 
of $K$~\cite{Tri02}. As pointed out by Nakamura~\cite{Nak04}, 
a possible approach to proving this for prime alternating links would 
be to establish $\max\deg_z P_{W(K)} = 2n-2$. 
In this special case, this statement is equivalent 
to Conjecture~\ref{conj:Kidwell}. We will study here a 
variant of this conjecture that also applies to links:
\begin{conj}\label{conj:extended}
Let $A(K,f)$ be a framed link. Then
$$\max\deg_z \{K,f\} = 2 \max\deg_x F_K$$
\end{conj}

We show that under an apparently weak additional condition 
on the knots under consideration,
the statement of Conjecture~\ref{conj:extended} implies 
the statement of Conjecture~\ref{conj:Kidwell}:

\begin{lem}\label{lem:bridge}
Let $A(K,f)$ be a framed knot for which  
$\max\deg_x F_K > 0$, and \\
$\max\deg_z \{K,f\} = 2 \max\deg_x F_K.$
Let $W_K$ be a Whitehead double of $K$.
Then Conjecture~\ref{conj:Kidwell} holds for $K$.
\end{lem}
\begin{proof}
Let $K$ be a knot with $\max\deg_x F_K > 0$ for which 
Conjecture~\ref{conj:extended} holds. Then for every framing $f$, 
we have $\max\deg_z [K,f] > 0$ since $\{K,f\} = [K,f]-1$. 
Let $W(K)$ be the $n$-twisted Whitehead double of $K$.
Applying a skein relation in the region of the clasp of the Whitehead double,
we get
$P_W(K) = v^2 \delta_P +vz [K,n]$, hence
$\max\deg_z P_{W(K)}  = 1 + \max\deg_z [K,n]$.
\end{proof}

The above additional condition $\max\deg_x F_K > 0$ appears to be 
very weak, since to the author's knowledge no example of a nontrivial knot 
with~$\max\deg_x F_K = 0$ is known. 
As a warm-up, it is easy to see that Conjecture~\ref{conj:extended}
in particular holds for all unlinks:
\begin{lem}\label{lem:unlink}
Conjecture~\ref{conj:extended} holds for all trivial links.
\end{lem}
\begin{proof}
Assume $f$ is the framing induced by a diagram $D$ for the $n$-component 
trivial link $\bigcirc^n$. Then
$$\{\bigcirc^n,f\} =v^{2sw(D)}\{\bigcirc\}^n =v^{2sw(D)}(1-\delta_P)^n$$
by elementary properties of the Rudolph polynomial stated in \cite{Rud90}.
Thus $\max\deg_z \{\bigcirc^n,f\}=0$. For the Kauffman polynomial,
we similarly have $F_{\bigcirc^n} = \delta_F^n$, and so the maximal 
$x$-degree of $F$ equals~$0$ as well.
\end{proof}

We briefly review some previous results related to the above conjectures.
The following generalizes an observation made by Tripp in~\cite{Tri02}:
\begin{prop}\label{prop:extended-ineq}
Let $K$ be a nontrivial $n$-crossing prime alternating link. 
Then for every framing $f$, the inequality
$\max\deg_z \{K,f\} \le 2 \max\deg_x F_K$ holds.
\end{prop}
\begin{proof}
Let $D$ be a reduced alternating diagram of $K$ with $n$ crossings, 
and let $\widetilde D$ be the diagram obtained by  
replacing each arc in $D$ by a corresponding pair of antiparallel arcs.

We apply Morton's inequality to $\widetilde D$: 
The number of crossings $c(\widetilde D)$ equals $4n$.
For counting the number of Seifert circles, we follow an idea given 
in~\cite{Tri02}: Regard $D$ as a $4$-regular planar graph with~$n$ 
vertices at the crossings, and, by Euler's polyhedral formula, 
with $n+2$ faces.
Then the number of Seifert circles $s(\widetilde D)$ arising from smoothing 
out the crossings equals $2n+2$: One for each face of the graph, and one for 
each four-crossing region in $\widetilde D$ that corresponds 
to a single crossing in $D$.
Since the maximum $z$-degree is not altered by changing the
framing, Morton's inequality gives $\max\deg_z [K,f] \le 2n-2$ for every 
framing $f$. Applying the same procedure after deleting some link 
components from the diagram~$D$ and corresponding components from the diagram
$D'$, we also get $\max\deg_z [K',f] < 2n-2$ for every proper sublink $K'$ of $K$.
Thus, $\max\deg_z \{K\} \le 2n-2$. On the other hand, Theorem~\ref{thm:Kidwell} 
implies that $\max\deg_x F_K = n-1$, since $D$ is a reduced prime alternating 
diagram.
\end{proof}
\begin{rem*}
A similar fact was claimed, but with an erroneous proof that 
is fixed here, in a draft of an earlier paper 
of the present author, see~\cite{Ich03}.
\end{rem*} 

Recently, Nakamura proved that $\max\deg_z [K,f] = 2n-2$ for every 
$n$-crossing rational link with framing $f$, 
generalizing a previous result by Tripp~\cite{Nak04,Tri02}. 
Since each of the proper sublinks of a $2$-component 
rational link $K$ is by itself a trivial knot or empty, we have in this case
$\{K\} = 1- 2\delta_P^2 + [K]$ and $\max\deg_z [K,f] =\max\deg_z \{K\} =
2n-2$. Together with Lemma~\ref{lem:bridge} and 
Proposition~\ref{prop:extended-ineq}, we obtain the following preliminary result:
\begin{prop}
Conjectures~\ref{conj:Kidwell} and~\ref{conj:extended} hold for all rational links.\quad\proofbox
\end{prop}
In the following, we will generalize this result to a much larger class 
of links.

\section{Split Union and Composition of Links.}

In this section, we prove that the relation from
Conjecture~\ref{conj:extended} is preserved under split union 
and composition, hence reducing the task to proving it 
for the case of prime links. The case of disjoint union is 
not hard to prove:

\begin{prop}\label{prop:split-links}
If Conjecture~\ref{conj:extended} holds for two framed links 
$L_1$, $L_2$, it also holds for the disjoint union $L_1\sqcup L_2$.
\end{prop}
\begin{proof}
For the Kauffman polynomial, we have
$F_{L_1\sqcup L_2}= F_{L_1}F_{L_2}$, see~\cite{Kau90};
and similar for the framed Rudolph polynomial 
$\{L_1\sqcup L_2,f\}=\{L_1,f\}\{L_2,f\}$, see~\cite{Rud90}.
\end{proof}

\begin{figure}\label{fig:numdensum}
\includegraphics[scale=.8]{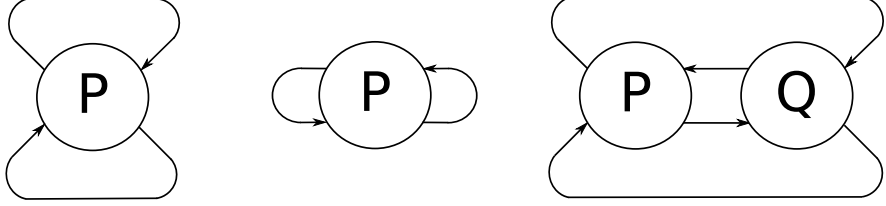}
\caption{Numerator closure, denominator closure and total sum}
\end{figure}

Next, we examine the case of composite links.
Conjecture~\ref{conj:Kidwell} suggests that the 
$z$-degree of $[L_1\#L_2,f]$ is additive for nontrivial knots
$L_1$,$L_2$ --- since for the Kauffman polynomial of composite links, 
the relation $F_{L_1\#L_2}=\delta_F F_{L_1}F_{L_2}$ holds~\cite{Kau90}.
But if one of the factor links $L_1$,$L_2$ is the unknot $\bigcirc$, 
the additive behavior of $\max\deg_z$
holds {\em if and only if} we choose the framing 
in a way such that $[\bigcirc,f]$ has nonnegative degree. 
Despite these seeming complications, we can indeed derive a formula for
$[L_1\#L_2,f]$ using linear skein theory.
We need the following definition.
\begin{defn}
For an appropriately oriented tangle~$P$, 
its {\em numerator closure}, denoted by~$P^N$, 
and its {\em denominator closure}, denoted by~$P^D$,
are the links obtained by  closing up the pending ends 
as shown in the left and middle diagram of 
Figure~\ref{fig:numdensum}. The right part of Figure~\ref{fig:numdensum} 
shows the link obtained as the {\em total sum} of two tangles $P$ and $Q$, 
denoted by~$(P+Q)^N$.
\end{defn}

\begin{prop}\label{composite}
Let $A(L,0)=K_A\#_{B}J$ be a composite framed link 
obtained by joining the component $A$ of $K$ 
with the component $B$ of $J$. Then 
\begin{equation*}
(\delta_P^2-1)[L]= \delta_P^2 [K\setminus A]\cdot
[J\setminus B] +  [K]\cdot[J]- [K\setminus A]\cdot[J] - [K]\cdot[J\setminus B]
\end{equation*}
\end{prop}
\begin{proof}
Lickorish and Millet  showed that the following formula holds 
for the skein polynomial of the total sum of two tangles 
$Q$ and $R$, referred
to as the numerator-denominator formula \cite{LicMil87}. 
$$(\delta_P^3-\delta_P)P_{(Q+R)^N} = \delta_P\left(P_{Q^N} P_{R^N} + P_{Q^D}P_{R^D}\right)-
\left(P_{Q^N}P_{R^D}+P_{Q^D}P_{R^N}\right),$$
 We can interpret $\partial A(L,0)$ as the total sum of tangles $Q$ and $R$,
with $\partial A(K,0)=Q^D$, and $\partial A(J,0)=R^D$. 
Then we get $$Q^N=\partial A(K\setminus A,0)\,\sqcup\,\bigcirc,~R^N = \partial A(J\setminus B,0)\,\sqcup\, \bigcirc.$$
Applying the numerator-denominator formula to this total sum
yields the result.
\end{proof}

We have obtained a somewhat complex formula describing the polynomial $[L]$
under link composition. Astonishingly, the Rudolph polynomial 
admits a rather simple and natural formula for link composition.

\begin{thm}[The Composition Formula]\label{thm:composition}
Let $L=K\#J$ be a composite link with factors $K$ and $J$. 
Then $$\{\bigcirc\}\{K\#J\}=\{K\}\{J\}.$$
\end{thm}

\begin{proof}
Assume the components $A$ and $B$ are joined to form the composite link
$L=K_A\#_{B}J$. Then we use $C$ to denote the resulting joint component.
We order the $r = 2^{|L|-1}$ sublinks of $L$ containing the 
joint component $C$, and denote them by $M_1,M_2\ldots M_{r}$. 
Then, by definition of the Rudolph polynomial,
\begin{equation}\label{drigg1}
\{\bigcirc\}\{L\}= (\delta_P^2-1)\left(\sum_{k=1}^{r}\left((-1)^{|L|-|M_k|}[M_k]\right)-\{L\setminus C\}\right).
\end{equation}
Note that the term $-\{L\setminus C\}$ on the right hand side of 
Equation~\eqref{drigg1} has negative sign 
since $(-1)^{|L|} = -(-1)^{|L\setminus C|}$.

Each of the sublinks $M_k$ is of the form $K'_A\#_B J'$ 
for some sublinks $K' \subseteq K$ and $J'\subseteq J$.
Proposition~\ref{composite} allows to compute the skein polynomial
of these links in terms of sublinks of $K \sqcup J$. In order to do so,
we associate with each $M_k$ the following four sublinks of $K \sqcup J$:
\begin{itemize}
\item the sublink $L_{1,k}=(K'\setminus A)\sqcup(J'\setminus B)$, containing
  neither $A$ nor $B$,
\item the sublink $L_{2,k}= K' \sqcup J'$, containing both $A$ and $B$,
\item the sublink $L_{3,k}=(K'\setminus A)\sqcup J'$, containing $B$ but not $A$,
and
\item the sublink $L_{4,k}=K'\sqcup(J'\setminus B)$, containing $A$ but not $B$. 
\end{itemize}

Note that the above association partitions the set of sublinks
of $K \sqcup J$ into four sets 
the sense that $\bigcup_{j=1}^4\left(\bigcup_{k=1}^r L_{j,k}\right)$
equals the set of all sublinks of $K \sqcup J$. Thus, 
\begin{equation}\label{drigg2}
\{K\sqcup J\} =  \sum_{j=1}^{4}\sum_{k=1}^{r}
 (-1)^{|K|+|J|-|L_{j,k}|}[L_{j,k}].
\end{equation}
Fix a sublink $M_k$ for the moment. By Proposition~\ref{composite} and 
the fact that the skein polynomial is multiplicative under split union, 
we obtain:
\begin{equation}
(\delta_P^2-1)[M_k] =
 \delta_P^2 [L_{1,k}]+[L_{2,k}]  -[L_{3,k}] -[L_{4,k}].
\end{equation}
Observe that $|L_{1,k}|$ and $|L_{2,k}|$ both have parity different
from $|M_k|$, while $|L_{3,k}|$ and $|L_{4,k}|$ have the same parity
as $|M_k|$. Thus, for instance $(-1)^{|L|-|M_k|} = (-1)^{|L|+1-|L_{1,k}|}$.
Hence, when multiplying both sides of the above equation with
$(-1)^{|L|-|M_k|}$, using analogous considerations, 
we obtain the following equation:
\begin{equation}
\begin{split}
(\delta_P^2-1)(-1)^{|L|-|M_k|}[M_k] =\delta_P^2& (-1)^{|L|+1-|L_{1,k}|} [L_{1,k}]\\
+&(-1)^{|L|+1-|L_{2,k}|} [L_{2,k}] \\
+&(-1)^{|L|+1-|L_{3,k}|}[L_{3,k}] \\
+&(-1)^{|L|+1-|L_{4,k}|}[L_{4,k}] \\
\end{split}
\end{equation}
By adding the term $(\delta_P^2-1) (-1)^{|L|-|L_{1,k}|}[L_{1,k}]$ 
to both sides of the equation, we get:
\begin{equation}\label{drigg3}
(\delta_P^2-1)\left((-1)^{|L|-|M_k|}[M_k] +(-1)^{|L|-|L_{1,k}|} [L_{1,k}]\right)
= \sum_{j=1}^4 (-1)^{|L|+1-|L_{j,k}|}[L_{j,k}]
\end{equation}
Now we treat $k$ again as variable, and
sum up left hand sides and right hand sides of
Equation~\eqref{drigg3} for $k=1$ up to $r$. By virtue of 
the equality $|L|+1 = |K|+|J|$, we can use 
Equation~\eqref{drigg2} to simplify the resulting right-hand side, and
get:
\begin{equation}\label{drigg4}
(\delta_P^2-1)\left(\sum_{k=1}^r
\left((-1)^{|L|-|M_k|}[M_k]\right) +\sum_{k=1}^r\left((-1)^{|L|-|L_{1,k}|} [L_{1,k}]\right)\right)
= \{K\sqcup J\}
\end{equation}
A basic property of the Rudolph polynomial is 
$\{K\sqcup J\} = \{K\}\{J\}$ (see~\cite{Rud90}), so the right-hand side of 
Equation~\eqref{drigg4} equals the right-hand side of the equation
in the statement of the theorem. 

To simplify the left-hand-side of Equation~\eqref{drigg4}, 
observe that the set $\{\,L_{1,k} \mid 1 \le k \le r\,\}$
equals the set of sublinks of the link $L \setminus C$.
Obviously, $|L|=|L\setminus C|+1$, hence  
\begin{equation}
\sum_{k=1}^r(-1)^{|L|-|L_{1,k}|} [L_{1,k}]= -\{L \setminus C\}.
\end{equation}
Therefore, the left hand side of Equation~\eqref{drigg4}
equals the right hand side of Equation~\eqref{drigg1}. 
This completes the proof of the theorem.
\end{proof}

The composition formula reduces the problem of proving
Conjecture~\ref{conj:extended} to the case of prime 
links. In particular, it allows us to drop the primality 
condition from Proposition~\ref{prop:extended-ineq}:
\begin{cor}
If $K$ is an alternating knot, then 
$$\max\deg_z \{K \} \leq 2\max\deg_x F_K.\mathproofbox$$
\end{cor}

With the aid of the composition formula, it is now  
also easy to verify that Conjecture~\ref{conj:Kidwell}
holds for infinitely many non-alternating knots, since that 
conjecture has been verified in particular 
for all prime (also non-alternating) 
knots of crossing number at most $11$~\cite{KidSto03}.
We can of course construct infinitely many different 
composite knots whose prime
factors are non-alternating and have at most $11$ crossings each.
In these cases, Lemma~\ref{lem:bridge} is applicable, so this
fact extends to Conjecture~\ref{conj:extended}.

We conclude this section with a further 
property of the Rudolph polynomial which can be derived
from the composition formula:
Note that the expression $K\#J$ is ambiguous in the case
where at least one of $K,J$ is a link---the ambiguity disappears 
only if we write $K_A\#_{B}J$ instead. Thus, similar to the case of  
the skein polynomial~\cite{LicMil87}, we have
\begin{cor}
There exist infinitely many pairs of nonequivalent links with the same
Rudolph polynomial.\quad\proofbox
\end{cor}

\section{Algebraic Alternating Links}

In this section, we apply the congruence theorem to 
identify a large class of alternating links for which
Conjecture~\ref{conj:extended} holds true. 
This class in particular contains all algebraic alternating links.

We need some additional notions first. 
The framed Kauffman polynomial $L_{K}(a,x)$ can be equivalently
seen as a single-variable Laurent polynomial in $x$ over $\mathbb Z[a^{\pm 1}]$.
Then the coefficient of the maximal power in $x$ is a polynomial in $a$,
denoted by $X_{K}(a)$. (Here we use the framed version because
below we make use of some results on that polynomial from~\cite{Thi88}.) 
Let $c(K)$ denote the minimal crossing number of~$K$.
Now the congruence theorem can be used to show that
Conjecture~\ref{conj:extended} holds in the case 
where $X_K(a)$ does not vanish when reduced modulo $2$:

\begin{thm}
Assume $K$ is an alternating link with $X_K(a)\not\equiv_2 0$.
Then Conjecture~\ref{conj:extended} holds for $K$.
\end{thm}

\begin{proof}
By virtue of Proposition~\ref{prop:extended-ineq}, it suffices 
to show the inequality $\max\deg_z \{K\} \ge 2\max\deg_x F_K(x,z)$.
Also, by virtue of Lemma~\ref{lem:unlink} and the composition formula 
(Theorem~\ref{thm:composition}), 
we may assume $K$ is nontrivial and prime.
Write $F_K(a,x)\mod 2$ for the Kauffman polynomial
with coefficients reduced modulo $2$. 
Since $X_K(a)$ does not vanish when reduced modulo $2$, 
$\max\deg_x L_K(a,x) \mod 2 =\max\deg_x F_K(a,x)$. Moreover,
Theorem~\ref{thm:Kidwell} yields $\max\deg_x L_K(a,x)=n-1$
for a nontrivial prime alternating framed link $K$ with minimal 
crossing number~$n$.
By the congruence theorem, we can deduce $\max\deg_z \{K,f\} \ge 2(n-1)$,
and since the maximum $x$-degree of the framed Kauffman polynomial 
is independent of the framing provided it is greater than $0$, 
which is the case since $n \ge 2$, we have
$\max\deg_z \{K\} \ge 2(n-1) = 2\max\deg_x F_K(x,z)$.
\end{proof}

The first alternating knots violating the condition $X_K(a)\not\equiv_2 0$
are the knots $8_{16}$ and $8_{17}$. These are also among the first 
alternating knots which are not algebraic. Indeed,
Thistlethwaite investigated the polynomial~$X_K(a)$
for prime alternating links, and found the equality 
\[ X_K(a) = \kappa  (a^{-1}+a)z^{n-1},\]
where $n$ denotes the number of crossings in a reduced alternating diagram, 
and $\kappa$ denotes the {\em chromatic invariant} of a basic Conway polyhedron
into which~$K$ can be inscribed, see~\cite{Biggs74,Con69,Thi88} for precise 
definitions and more background information.
There are probably infinitely many basic Conway polyhedra with odd chromatic
invariant, and each of them gives rise to an infinite family of alternating 
links for which the above theorem is applicable. 
In particular, algebraic alternating links can be inscribed into the 
basic Conway polyhedron $1^*$, for which $\kappa=1$ holds~\cite{Thi88}. 
We thus obtain our second main result:
\begin{thm}
Conjecture~\ref{conj:extended} holds for all algebraic alternating links.
\quad\proofbox
\end{thm}
Furthermore, the same fact holds for alternating links that can be 
inscribed into the basic Conway polyhedra $8^*$, $10^{**}$, and a 
few other small basic polyhedra which are known to have odd 
chromatic invariant, as computed in~\cite{Thi88}.

\section{Conclusion}

To the author's knowledge, the relation between the skein polynomial
and the $2$-variable Kauffman polynomial is not yet understood very well.
A notable exception is Rudolph's Congruence Theorem. The present paper 
contains a framed variant of this theorem, which relates the 
framed Kauffman polynomial of a knotted annulus $A(K,f)$ and the 
skein polynomial of its boundary link $\partial A(K,f)$.

Next, a new conjectured relation between the skein polynomial
and the $2$-variable Kauffman polynomial was presented, 
which generalizes, at least under a weak additional assumption, 
a recent conjecture by Kidwell and Stoimenow, and, in the case of 
prime alternating knots, a recent conjecture by Tripp. 
The present research report studied the the behavior of 
the Rudolph polynomial under basic operations on links. In particular,
an unexpectedly natural formula for the Rudolph polynomial of composite links
was found. We proved that the conjectured relation is preserved under split 
union and link composition, and established its truth for a large class of 
alternating links. This class properly contains the algebraic 
alternating links, and for this class
the mentioned additional assumption is true, namely that the $x$-degree of the
Kauffman polynomial is greater than zero.
Hence, the conjectures by Kidwell and Stoimenow and Tripp turn out to hold
for algebraic alternating knots. A weaker result heading into the same direction
is found in a recent preprint by Brittenham and Jensen, establishing 
the conjecture for alternating pretzel knots~\cite{BriJen06}. 
That work, although appearing later than a preprint of the present paper,
is apparently independent of the present work,
since their proof is in part geometric, based on canonical Seifert surfaces.
It would be interesting to see whether these approaches admit a joint
generalization, in particular for proving the conjecture for a larger family 
of non-alternating links. 
Another interesting line of research concerns the relation of 
the framing variable~$v$ of the Kauffman polynomial and the 
Rudolph polynomial. 
Some results in this direction are found already in~\cite{Nutt97,Nutt99,Rud90}.

\begin{acknowledgment}
Many thanks to Lee Rudolph for pointing the author to  
reference~\cite{Rud90}; this enabled the author to establish 
substantially more general results than in a first draft of 
the present paper.
Thanks also to Alexander Stoimenow for some interesting discussion
during the early stage of this work, and to an anonymous referee for some
valuable corrections and suggestions.
\end{acknowledgment}

\bibliographystyle{numeric}

\begin{thebibliography}{9}

\bibitem{Biggs74}
{N. Biggs}. {\em Algebraic Graph Theory.}
Cambridge University Press, London, UK, 1974

\bibitem{BriJen06}
{M. Brittenham and J. Jensen}.
Canonical genus and the {Whitehead} doubles of pretzel knots.
available online at \url{http://arxiv.org} as \url{arXiv:math/0608765v1}, 
2006, 16 pages.

\bibitem{BZ85}
{G. Burde and H. Zieschang}. {\em Knots.}
Walter de Gruyter \& Co., Berlin, Germany, 1985.

\bibitem{Con69}
{J. H. Conway}. An enumeration of knots and links, and some of their
algebraic properties. In: D. D. Leech (ed.), {\em Computational Problems in
Abstract Algebra.} Pergamon Press, Oxford, UK, 1969, 329--358.

\bibitem{HOMFLY}
{P. Freyd, J. Hoste, W.~B.~R. Lickorish, K. Millett, A. Ocneanu and 
D. Yetter}. A new polynomial invariant of knots and links. 
{\em Bull. Amer. Math. Soc. (N.S.)} \textbf{12(2)} (1985), 239--246.  

\bibitem{Ich03}
{H. Gruber}. Estimates for the minimal crossing number.
available online at \url{http://arxiv.org} as \url{arXiv:math/0303273v3},
2003, 11 pages.

\bibitem{Jon85}
{V. F. R. Jones}. A polynomial invariant of knots and links via von {N}eumann
algebras. {\em Bull. Amer. Math. Soc. (N.S.)} \textbf{12(1)} (1985), 103--111.

\bibitem{Kau90}
{L. H. Kauffman}. An invariant of Regular Isotopy. {\em
  Trans. Amer. Math. Soc.} \textbf{318(2)} (1990), 417--471.

\bibitem{Kid87}
{M. E. Kidwell}. On the degree of the Brandt-Lickorish-Millett-Ho polynomial
of a link. 
{\em Proc. Amer. Math. Soc.} \textbf{100} (1987), 755--762.

\bibitem{KidSto03}
{M. Kidwell and A. Stoimenow}. Examples related to the crossing number,
writhe, and maximal bridge length of knot diagrams. 
{\em Michigan Math. J.} \textbf{51(1)} (2003), 3--12.

\bibitem{Lic88}
{W. B. R. Lickorish}. Polynomials for Links. {\em Bull. London Math. Soc.}
\textbf{20(6)} (1988), 558--588.

\bibitem{LicMil87}
{W. B. R. Lickorish and K. Millett}. A polynomial invariant of oriented
links. {\em Topology} \textbf{26(1)}, 1987, 107--141.

\bibitem{Mor86}
{H. R. Morton}. {S}eifert circles and knot polynomials. 
{\em Math. Proc. Camb. Phil. Soc.} \textbf{99}, 1986, 107--109.

\bibitem{Nak04}
{T. Nakamura}. On the crossing number of $2$-bridge knot and the canonical
genus of its {W}hitehead double. 
{\em Osaka J. Math.} \textbf{43(3)} (2004), 609--623.

\bibitem{Nutt97}
{I. Nutt}. Arc index and the {Kauffman} polynomial. 
{\em J. Knot Theory and its Ramifications} \textbf{6(1)} (1997), 61--77.

\bibitem{Nutt99}
{I. Nutt}. Embedding knots and links in an open book {III}: On the braid index
of satellite links. {\em Math. Proc. Camb. Phil. Soc.} \textbf{126} (1999), 77--98.

\bibitem{Oht04}
{T. Ohtsuki (ed.)}.
Problems on invariants of knots and 3-manifolds. In:
T. Ohtsuki, T. Kohno, T. Le, J. Murakami, J. Roberts and V. Turaev (eds.),
{\em Invariants of knots and 3-manifolds (Kyoto 2001).} 
Geometry \& Topology Monographs \textbf{4} (2004), 377--572.

\bibitem{PT}
{J. Przytycki and P. Traczyk}. Conway Algebras and Skein Equivalence of Links.
{\em Proc. Amer. Math. Soc.} \textbf{100} (1987), 744-748.

\bibitem{Rud90}
{L. Rudolph}. A congruence between link polynomials.
{\em Math. Proc. Camb. Phil. Soc.} \textbf{107} (1990), 319--327.

\bibitem{Thi88}
{M. Thistlethwaite}. {K}auffman's polynomial and alternating links.
{\em Topology} \textbf{27(3)} (1988), 311--318.

\bibitem{Tri02}
{J. Tripp}. The canonical genus of an infinite family of knots.
{\em J. Knot Theory and its Ramifications} \textbf{11(8)} (2002),1233--1242.

\bibitem{Yam89}
{S. Yamada}. An operator on regular isotopy invariants of link diagrams.
{\em Topology} \textbf{28(3)} (1989), 369--377.

\end{thebibliography}

\end{document}